\documentclass[a4paper]{amsart}
\usepackage[foot]{amsaddr}
\usepackage{mathtools}
\usepackage{url}
\providecommand\given{}
\newcommand\SetSymbol[1][]{%
  \nonscript\:#1\vert
  \allowbreak
  \nonscript\:
  \mathopen{}}
\DeclarePairedDelimiterX\Set[1]{\{}{\}}{%
  \renewcommand\given{\SetSymbol[\delimsize]}
  #1}

\DeclarePairedDelimiterXPP\pospart[1]{}{(}{)}{^+}{#1}
\DeclarePairedDelimiterXPP\negpart[1]{}{(}{)}{^-}{#1}
\DeclarePairedDelimiter\ceil{\lceil}{\rceil}

\DeclarePairedDelimiter\interv{[}{]}
\DeclarePairedDelimiter\openinterv{(}{)}
\DeclarePairedDelimiter\absval{\lvert}{\rvert}

\newcommand\R{\mathbb{R}}
\newcommand\Z{\mathbb{Z}}
\newcommand\N{\mathbb{N}}
\newcommand\D{\mathcal{D}}

\newcommand\mapping[3]{{#1}\colon{#2}\to{#3}}
\newcommand\transp[1]{{#1}^{\scriptscriptstyle \mathsf{T}}}

\newtheorem*{prop0}{Proposition}

\newtheorem*{thm0}{Theorem}

\DeclareMathOperator*{\Int}{int}

\DeclareMathOperator*{\conv}{co}

\makeatletter
\newcommand{\leqnomode}{\tagsleft@true\let\veqno\@@leqno}
\newcommand{\reqnomode}{\tagsleft@false\let\veqno\@@eqno}
\makeatother

\usepackage{xifthen}
\newcommand\datapoints{X}
\newcommand\datapoint{\lowercase\expandafter{\datapoints}}
\newcommand\data[1][]{     %
  \datapoint\ifthenelse{\isempty{#1}}{}{^{#1}} %
}
\newcommand\ndata{N}
\newcommand\ordcone{C}
\newcommand\valf[2][]{\phi_{\ifthenelse{\isempty{#1}}{\datapoints}{{#1}},\prob}
  \ifthenelse{\isempty{#2}}{}{({#2})} %
}
\newcommand\valfAlt[2][]{\tilde{\phi}_{\ifthenelse{\isempty{#1}}{\datapoints}{{#1}},\prob}
  \ifthenelse{\isempty{#2}}{}{({#2})} %
}
\newcommand\prob{p}
\newcommand\quantile[2][]{     %
  q^{-}_{\ifthenelse{\isempty{#1}}{\datapoints}{#1}} %
  \ifthenelse{\isempty{#2}}{}{(#2)}
}
\newcommand\cquantile[2][]{     %
  Q^{-}_{\ifthenelse{\isempty{#1}}{\datapoints,\ordcone}{#1}} %
  \ifthenelse{\isempty{#2}}{}{(#2)}
}

\title[Finite Representation of Quantile Sets]{Finite Representation
  of Quantile Sets for Multivariate Data via Vector Linear
  Programming}

\author{Andreas L{\"o}hne}
\address[Andreas Löhne]{Friedrich Schiller University Jena, Germany}
\email{andreas.loehne@uni-jena.de}

\author{Benjamin Weißing}
\address[Benjamin Weißing]{Free University of Bozen--Bolzano, Italy}
\email{benjamin.weissing@unibz.it}

\begin{document}


\begin{abstract} Empirical quantiles for finitely distributed
  univariate random variables can be obtained by solving a certain
  linear program. It is shown in this short note that multivariate
  empirical quantiles can be obtained in a very similar way by solving
  a vector linear program. This connection provides a new approach for
  computing Tukey depth regions and more general cone quantile sets.

  \smallskip
  \noindent
  \textbf{\keywordsname:} multivariate quantiles, vector linear
  programming, multiple objective linear programming
\end{abstract}
\maketitle
\noindent
Let
\(\datapoints = \Set*{\data[1], \data[2], \ldots, \data[\ndata]}
\subseteq \R\) be a finite set of data points. The
\emph{empirical lower quantile function} of \(\datapoints\) is given by
\reqnomode
\begin{equation}
  \label{eq:scquantile}
  \mapping{\quantile[]{}}{\openinterv{0,1}}{\R}
  \text{,}\qquad
  \quantile[]{\prob}
  \coloneq
  \min\Set*{
    \bar{\datapoint} \in \datapoints
    \given
    \#\Set{
      \datapoint \in \datapoints
      \given
      \datapoint \leq \bar{\datapoint}
    } \geq \ceil{\ndata \prob}
  }\text{,}
\end{equation}
where \(\ceil{y}\coloneq \min\Set*{z \in \Z \given z\geq y}\) denotes
the smallest integer which is not smaller than \(y\).\par
Instead of using \eqref{eq:scquantile} directly,
\(\quantile[]{\prob}\) can also be computed by solving the
minimization problem
\begin{equation}
  \label{eq:lp_d_alt}
  \min\valfAlt{t} \qquad \text{s.t.} \quad t \in \R\text{,}
\end{equation}
where \(\mapping{\valfAlt[]{}}{\R}{\R}\) is defined as
\begin{equation*}
  \label{eq:phi}
  \valfAlt[]{t} \coloneq
  \sum_{i = 1}^\ndata
  \prob\pospart{\data[i] - t}
  +
  (1 - \prob)\negpart{\data[i] - t}\text{,}
\end{equation*}
with \(y^{+}\coloneq\max\Set{0, y}\) and \(y^{-}\coloneq\max\Set{0,
  -y}\) for \(y\in\R\), cf.\ \cite{KoeBas78}.
Under the mild assumption of non-integrality of \(\ndata\prob\), the
\(\prob\)-quantile of \(\datapoints\) is equal to the optimal solution
of \eqref{eq:lp_d_alt}, see \cite[Theorem~3.3]{KoeBas78} and the
following remark therein.
\begin{prop0}\label{prop:1}
  If \(\ndata p \notin \Z\), then \eqref{eq:lp_d_alt} has the
  unique optimal solution \(\bar{t} = \quantile[]{\prob}\).
\end{prop0}
The proof to this result can be found in \cite{KoeBas78}.  We consider
a slightly altered variant of the objective function:
\begin{equation*}
  \valf{t} \coloneq \valfAlt{t} + (1 - \prob)\sum_{i = 1}^\ndata\data[i]\text{.}
\end{equation*}
While the addition of a constant does not alter the optimal solution,
it simplifies the description of the dual problem discussed below
considerably.  Hence the departure point of our exposition is the
optimization problem
\leqnomode
\begin{equation}
  \label{eq:lp_d}
  \tag{LP\ensuremath{^*}\kern-\scriptspace}
  \min\valf{t} \qquad \text{s.t.} \quad t \in \R\text{.}
\end{equation}
\reqnomode
The objective function \(\valf[]{}\) of \eqref{eq:lp_d} is
\emph{polyhedral convex} (e.g., piecewise affine), meaning its
epigraph \(\Set{\openinterv{t,\, r} \in \R^2 \given \valf{t}\leq r}\) is a convex
polyhedron.  By introducing auxilliary variables, \eqref{eq:lp_d} may
be transformed into a linear program.  A solution to \eqref{eq:lp_d}
can thus be obtained by solving the corresponding dual (compare
\cite[Appendix]{KoeBas78}):
\leqnomode
\begin{equation}
  \label{eq:lp}\tag{LP}
  \max \sum_{i=1}^\ndata  \data[i]u_i \qquad \text{s.t.}\quad
  \begin{cases}
    \sum\limits_{i = 1}^\ndata u_i = (1 - \prob) \ndata\\
    0 \leq u_i\leq 1 & i \in \Set{1, \ldots, \ndata}\\
    u \in \R^\ndata
  \end{cases}
  \text{.}
\end{equation}
\reqnomode\par
In this short note, we show that this relationship between univariate
quantiles and linear programming can be extended to
\emph{multivariate} cone quantiles and \emph{vector} linear
programming.
Multivariate cone quantiles are introduced in \cite{HamKos18} and
extend the notion of quantiles to multivariate variables
\(\datapoints\subseteq\R^d\).  The well developed framework of vector
linear programs (VLP) extends scalar linear programming to
vector-valued objective functions and provides a duality theory called
\emph{geometric duality} (see \cite{HeyLoe08} and compare
\cite{LoeWei17}).  An important aspect of geometric duality in the
current context is its equivalence to classical LP-duality.\par
While applying the definition from \eqref{eq:scquantile} to compute
the \(\prob\)-quantile of a univariate data~set is straightforward
(sort the elements \(\data[i]\) in ascending order, then select the
\(\ceil{\ndata\prob}\)-th element of the sorted data~set), an imminent
difficulty in carrying over this procedure to multivariate quantiles
lies in the absence of a total ordering in \(\R^d\) and the ensuing
inability of sorting data~sets consisting of vectors.  On the other
hand it is quite easy to find a multivariate analogy of \eqref{eq:lp}:
Replacing the univariate \(\data[i]\) in \eqref{eq:lp} by vectors
\(\data[i]\in\R^d\) directly results in a VLP.  We employ geometric
duality for VLPs in order to retrieve a vector-valued counterpart of
\eqref{eq:lp_d}.  We show that the solution to this dual VLP can,
analogous to the univariate case, be used to construct the
multivariate \(\prob\)-quantile.  The major contribution of this
article is the ability to compute multivariate quantiles by using a
VLP-solver like the one presented in \cite{LoeWei17}.\par
For both concepts, VLPs as well as multivariate cone quantiles, the
notion of \emph{order cone}s is crucial.  Although vectors in \(\R^d\)
may not be sorted in general, it is possible to sort them at least
\emph{partially} by equipping \(\R^d\) with a \emph{partial order}.
To this end a polyhedral convex cone \(C \subseteq \R^d\), which does
not contain lines, is considered.  The cone \(C\) defines a partial
order (reflexive, transitive and antisymmetric relation) “\(\leq_C\)”
via
\begin{equation}\label{eq:pord}
  y^1 \leq_C y^2 \ratio\iff y^2 - y^1 \in C\text{.}
\end{equation}
With this notion, one can define \emph{\(C\)-minimal} elements of a
set \(A\subseteq\R^d\): An element \(\bar{y}\in A\) is \(C\)-minimal
iff there is no element in \(A\) which is smaller:
\begin{equation*}
  \interv{y\in A\colon y \leq_C \bar{y}} \implies y = \bar{y}\text{.}
\end{equation*}
\emph{\(C\)-maximal} elements \(\bar{y} \in A\) are defined analogously.
The \emph{positive dual cone} of \(C\) is the cone
\(C^{+}\coloneq\Set{w \in \R^d \given y \in C \implies \transp{w}y
  \geq 0}\).  In the following, we assume that \(C\) has non-empty
interior.  Let \(c\) be an interior point of \(C\). Then
\(\transp{w}c > 0\) for all \(w \in C^{+}\) and
\begin{equation}
  \label{eq:cbasis}
  B^{+} \coloneq \Set{w \in C^{+} \given \transp{w}c = 1}
\end{equation}
defines a \emph{basis} of \(C^{+}\).  Observe that \(B^{+}\) is a
\(\openinterv{d - 1}\)-dimensional polytope and an element
\(w\in B^{+}\) is uniquely determined by its first \(d - 1\)
components.\par
The general idea of quantiles for multivariate data
\(\datapoints = \Set{\data[1], \data[2], \ldots, \data[\ndata]}
\subseteq \R^d\) is to apply formula~\eqref{eq:scquantile} to the
\emph{scalarized} data~set
\(\transp{w}\datapoints \coloneq \Set*{\transp{w}\data[1],
  \transp{w}\data[2], \ldots, \transp{w}\data[\ndata]} \subseteq \R\)
for some vector \(w \in B^{+}\).  Using \(w \in B^{+}\) as
\emph{weight vectors} preserves the partial order \eqref{eq:pord}
generated by \(C\): for all \(w \in B^+\),
\(\data[i] \leq_C \data[j]\) implies
\(\transp{w}\data[i] \leq \transp{w}\data[j]\).
For \(d \in \N\) and \(p \in \openinterv{0,1}\), the \emph{empirical
  lower \(C\)-quantile} of \(X\), introduced in a more general form in
\cite{HamKos18}, is the set
\begin{equation}
  \label{eq:2}
  \cquantile[]{\prob}
  \coloneq
  \bigcap_{w \in B^{+}}\Set*{
    z \in \R^d
    \given
    \transp{w} z \geq \quantile[\transp{w}\datapoints]{\prob}
  }\text{.}
\end{equation}
For simplicity, we use here the characterization from
\cite[Proposition 6]{HamKos18} as a definition. For other equivalent
variants, the reader is referred to \cite{HamKos18}. For the empirical
variants, see also \cite{HamKos22}.

In addition to the parameters \(\datapoints\) and \(\prob\) used in
\eqref{eq:scquantile}, a third parameter occurs in \eqref{eq:2}: The
lower \(C\)-quantile depends on the order cone \(C \subseteq \R^d\).
In the univariate setting (\(d = 1\)), this parameter is not necessary
as there exist essentially (excluding the trivial cones \(\Set{0}\)
and \(\R\)) two order cones in \(\R\): The non-negative and
non-positive real numbers \(C = \R_{+}\) (corresponds to \({\leq}\))
and \(C = \R_{-}\) (corresponds to \({\geq}\)), respectively.

As we will also see by our main result, in the intersection in
\eqref{eq:2} finitely many \(w \in B^{+}\) are known to be sufficient
(see also \cite{HamKos22}) for generating \(\cquantile[]{\prob}\).
The question is how to determine these weights \(w\).  For the
bivariate case (\(d = 2\)) an algorithm for computing
\(\cquantile[]{\prob}\) is stated in \cite{HamKos22}.  While the
algorithm provided there is not readily extendable to dimensions
\(d > 2\), the approach we use here is valid for arbitrary dimension
\(d\).

For a given multivariate data~set
\(\datapoints = \Set*{\data[1], \data[2], \ldots, \data[\ndata]}
\subseteq \R^d\), an order cone \(C \subseteq \R^d\) and
\(\prob\in\openinterv{0, 1}\), we consider the following VLP:
\leqnomode
\begin{equation}\label{eq:vlp}\tag{VLP}
  \max_{C} \sum_{i=1}^\ndata  \data[i]u_i \qquad \text{s.t.}\quad
  \begin{cases}
    \sum\limits_{i = 1}^\ndata u_i = (1 - \prob) \ndata\\
    0 \leq u_i\leq 1 & i \in \Set{1, \ldots, \ndata}\\
    u \in \R^\ndata
  \end{cases}\text{.}
\end{equation}
\reqnomode
Note that in contrast to \eqref{eq:lp} the objective function in
\eqref{eq:vlp} is vector-valued while the constraint sets of both
problems coincide.  A solution to \eqref{eq:vlp} is a finite \emph{set} of
feasible vectors \(u\in\R^\ndata\) rather than a single point.  The
order cone \(C\) in \eqref{eq:vlp} indicates that those elements \(u\)
of a solution are required to be \emph{\(C\)-maximizers}.  A
\(C\)-maximizer is a feasible point whose objective value is
\(C\)-maximal (compare \eqref{eq:pord}) among the set of all possible
objective values.  Hence, when \(d = 1\) and \(C = \R_{+}\), both
problems \eqref{eq:vlp} and \eqref{eq:lp} are equivalent.

As \eqref{eq:vlp} is used merely as a means to an end, namely for
constructing its dual, further aspects of the solution concept for
VLPs will be discussed in the context of this dual problem.  The
interested reader may refer to \cite{LoeWei17} and the references
therein for a comprehensive exposition of the VLP solution concept.

Constructing the dual of \eqref{eq:vlp} requires the selection of
\(c \in \Int C\), the so-called \emph{duality parameter}.  Without
loss of generality, we can choose \(c\) such that the last component
\(c_d\) is non-zero, i.e.\ \(\absval{c_d}\neq 0\). (At least one of
the components, say \(c_k\), is different from zero, as \(0\) is not
an interior point of \(C\).  In case \(c_d = 0\), we can exchange the
coordinate axes \(k\) and \(d\).)  Recall that \(c\) defines a basis
\(B^{+}\) of \(C^{+}\) by means of \eqref{eq:cbasis}.  The
\emph{geometric dual} of \eqref{eq:vlp} is
\begin{equation}\label{eq:geom_dual}
	\min_K \begin{pmatrix}
    \tfrac{c_d}{\absval{c_d}} w_1 \\
    \vdots \\
		\tfrac{c_d}{\absval{c_d}} w_{d-1} \\
    t \ndata (1 - \prob) +
    \sum\limits_{i = 1}^\ndata r_i
  \end{pmatrix} \; \text{ s.t. } \,
  \begin{cases}
    w \in B^+ &\\
    r_i \geq \transp{w}\data[i] - t &\text{for \(i \in \Set{1,\dots,\ndata }\)}\\
    r_i \geq 0 &\text{for \(i \in \Set{1,\dots,\ndata }\)}\\
    \mathrlap{\left(r,\,w,\,t\right) \in \R^\ndata \times \R^d \times \R}
  \end{cases}\text{,}
\end{equation}
with order cone
\begin{equation*}
  K \coloneq \Set{y \in \R^d \given y_1 = \cdots = y_{d - 1} =
    0,\, y_d \geq 0}.
\end{equation*}

Geometric duality is introduced in \cite{HeyLoe08}. The slightly
extended variant we use here can be found in
\cite[(VLP\textsubscript{max}) (note the typo that in the dual problem
\(K\)-maximize has to be replaced by \(K\)-minimize)]{LoeWei17}.
Problem \eqref{eq:geom_dual} is a vector linear program again; in
contrast to \eqref{eq:vlp} we are interested in \(K\)-\emph{minimal}
points and the corresponding feasible points generating them,
\(K\)-minimizers.  The dual cone \(K\) has a special property: Two
points \(y^1, y^2 \in \R^d\) are comparable under the order generated
by \(K\) (see \eqref{eq:pord}) if and only if they coincide in their
first \(d - 1\) components.  This means that to every \(K\)-minimal
point \(y\) of \eqref{eq:geom_dual} there exists \(\bar{w}\in C^{+}\)
such that
\begin{align}
  y_i &= \frac{c_d}{\absval{c_d}} \bar{w}_i
        \qquad
        \text{for \(i\in\Set*{1,\ldots, d - 1}\),}\nonumber
        \intertext{and}
        y_d &= \min\Set*{
              t \ndata (1 - \prob) + \sum\limits_{i = 1}^\ndata r_i
              \given
              (r, t)\in \R^\ndata \times \R\colon (\bar{w}, r, t)\in T
              }\text{,}\label{eq:minlastcomp}
\end{align}
where we denote by \(T\) the feasible set of \eqref{eq:geom_dual},
hold.  The two constraints in \eqref{eq:geom_dual} involving \(r_i\)
can be conflated to \(r_i \geq \pospart{\transp{w}\data[i] - t}\) for
\(i \in \Set{1, \ldots, \ndata}\).  In order to attain the minimum in
\eqref{eq:minlastcomp}, actually
\(r_i = \pospart{\transp{\bar{w}}\data[i] - t}\) needs to hold.  With
this in mind and taking into account
\begin{align*}
  t\ndata (1 - \prob) +
  \sum\limits_{i = 1}^\ndata\pospart{\transp{w}\data[i] - t}
  &=
  \valfAlt[\transp{w}\datapoints]{t} +
    (1 - \prob)\sum\limits_{i = 1}^\ndata \transp{w}\data[i]\\
  &=
    \valf[\transp{w}\datapoints]{t}
  \text{,}
\end{align*}
(where we use the polyhedral convex function
\(\valf[\transp{w}\datapoints]{}\) as introduced above in problem
\eqref{eq:lp_d}) the dual problem \eqref{eq:geom_dual} can be restated
more compactly as
\leqnomode
\begin{gather}\label{eq:vlp_d}\tag{VLP\ensuremath{^*}\kern-\scriptspace}
  \min_K  D(w, t ) \qquad \text{s.t.\ \((w,\,t) \in B^{+}\times \R\)}
\end{gather}
\reqnomode
with dual objective function
\begin{equation*}
  D(w, t) \coloneq
  \begin{pmatrix}
  	\tfrac{c_d}{\absval{c_d}} w_1 \\
	  \vdots \\
    \tfrac{c_d}{\absval{c_d}} w_{d-1} \\
    \valf[\transp{w}\datapoints]{t}
  \end{pmatrix}
  \text{.}
\end{equation*}
Now denote by \(T \coloneq B^{+} \times \R\) the feasible set of
\eqref{eq:vlp_d} and let
\(D[T] \coloneq \Set{D(w,\, t) \given (w,\, t) \in T}\) be its image.
Then the \emph{extended image} \(\D \coloneq D[T] + K\) is a convex
polyhedron.  A \emph{solution} to \eqref{eq:vlp_d} is a finite set
\(\bar{T}\subseteq T\) of \emph{\(K\)-minimizers} such that the
convex hull of \(D[\bar{T}]\) generates the whole extended image:
\begin{equation}
  \label{eq:convhulld}
  \D = \conv D[\bar{T}] + K\text{.}
\end{equation}
Equation \eqref{eq:convhulld} may be interpreted as a so-called
vertex-representation of the convex polyhedron \(\D\).  Because all
vertices are included in such a representation, we get the following
important implication of \eqref{eq:convhulld}: for every vertex \(y\)
of \(\D\) there exists a minimizer \((w, t) \in \bar{T}\) from the
solution set such that \(y = D(w, t)\).\par
Another useful observation concerns the structure of \(K\)-minimizers
of \eqref{eq:vlp_d}: for a \(K\)-minimal point \(y\) of
\eqref{eq:vlp_d}, we have (compare \eqref{eq:minlastcomp})
\begin{equation}\label{eq:kmin}
  y_d = \min\Set{\valf[\transp{w}\datapoints]{t}\given t\in \R}\text{.}
\end{equation}
As this minimum is attained uniquely in
\(\quantile[\transp{w}\datapoints]{\prob}\), any minimizer of
\eqref{eq:vlp_d} has the form
\((w, \quantile[\transp{w}\datapoints]{\prob})\) for some element
\(w \in B^{+}\).\par
With these preparations we are now ready to state the main result.  As
in the univariate case, a solution to \eqref{eq:vlp_d} can be used to
compute the \(\prob\)-quantile \(\cquantile[]{\prob}\).
\begin{thm0} Let
  \(X = \Set*{\data[1],\data[2],\ldots,\data[\ndata]} \subseteq
  \R^d\), \(\prob \in \openinterv{0, 1}\) such that
  \(\prob\ndata \notin \Z\), \(C \subseteq \R^d\) a line-free
  polyhedral convex cone with interior point \(c\) such that
  \(c_d \neq 0\). Then \eqref{eq:vlp_d} has a solution \(\bar {T}\)
  which provides the finite representation
  \begin{equation}\label{eq:main}
    \cquantile[]{\prob}
    =
    \bigcap_{(w,t) \in \bar{T}}
    \Set*{z \in \R^d
      \given
      \transp{w} z
      \geq t
    }\text{.}
  \end{equation}
\end{thm0}
\begin{proof}
  As \eqref{eq:vlp} has a non-empty, bounded feasible set, a solution
  to \eqref{eq:vlp} exists by \cite[Corollary~6]{LoeWei16}.  The
  existence of a solution to \eqref{eq:vlp_d} then follows by
  geometric duality \cite{HeyLoe08}.\par
  The first inclusion \(\subseteq\) follows directly from the
  observation made earlier: every element \((w, t) \in \bar{T}\) is of
  the form \((w, \quantile[\transp{w}\datapoints]{\prob})\), with
  \(w \in B^{+}\).\par
  For the reverse inclusion \(\supseteq\), we assume that the right
  hand side of \eqref{eq:main} is non-empty.  (Otherwise the inclusion
  is trivial.)  Let $z \in \R^d$ such that
  \begin{equation}
    \label{eq:pr0}
    \transp{w} z \geq t \qquad\text{holds
    for all \((w,t)\in \bar T\).}
  \end{equation}
  According to the definition of cone quantiles \eqref{eq:2} we need
  to show that
  \(\transp{w} z \geq \quantile[\transp{w}\datapoints]{\prob}\) holds
  for all \(w \in B^{+}\).  Take \(\bar{w} \in B^+\).  As noted above,
  the \(K\)-minimal point \(\bar{y} \coloneq D(\bar{w}, \bar{t})\)
  with
  \(\bar{t} \coloneq \quantile[\transp{\bar{w}}\datapoints]{\prob}\)
  corresponds to \(\bar{w}\).  As \(\bar{y}\) is \(K\)-minimal, it
  belongs to the relative interior of a \(K\)-minimal face of \(\D\).
  We denote by \(y^1, \ldots, y^k\) the vertices of this face, which
  are vertices of \(\D\) in particular and admit the representation
  \(y^j = D(w^j, t^j)\), with certain \((w^j,\, t^j) \in \bar{T}\),
  \(w^j\in B^{+}\) and
  \(t^j \coloneq \quantile[\transp{w^j}\datapoints]{\prob}\) (compare
  \eqref{eq:convhulld}).  The point \(\bar{y}\) is a convex
  combination of these vertices, i.e.,\ there exist
  \(\lambda_1, \ldots, \lambda_k \in \R\) with \(0\leq \lambda_j\) and
  \(\sum_j\lambda_j = 1\) such that
  \begin{equation}
    \label{eq:convcomb}
    \bar{y} = \sum_{j = 1}^k \lambda_j y^j\text{.}
  \end{equation}
  We now show that the preimage of \(\bar{y}\),
  \(\openinterv{\bar{w},\,\bar{t}}\), is a convex combination of the
  preimages of \(y^j\), \(\openinterv{w^j,\,t^j}\), with the same
  coefficients as in \eqref{eq:convcomb}.  For the first \(d - 1\)
  components of \(\bar{w}\) this follows directly from the definition
  of the dual objective function in \eqref{eq:vlp_d}.  For the
  \(d\)-th component the equality follows with \eqref{eq:cbasis},
  hence \(\bar{w} = \sum_{j = 1}^k \lambda_j w^j\).
  Finally we have to show \(\bar{t} = \sum_{j = 1}^k \lambda_j t^j\).
  To this end, we consider the extended VLP-formulation
  \eqref{eq:geom_dual} and construct the feasible point
  \((\bar{w},\, \sum_{j = 1}^k \lambda_j r^j,\, \sum_{j =
    1}^k\lambda_j t^j)\), where
  \(r^j_i \coloneq \pospart{\transp{w^j} \data[i] - t^j}\).  The
  \(d\)-th component of the corresponding image point is
  \begin{align*}
    \biggl(\,\sum_{j = 1}^k \lambda_j t^j\,\biggr) \ndata (1 - \prob) +
    \sum_{i = 1}^\ndata \sum_{j = 1}^k \lambda_j r^j_i
    &= \sum_{j = 1}^k \lambda_j \left(t^j \ndata (1 - p) + \sum_{i = 1}^\ndata r^j_i\right)\\
    &= \sum_{j = 1}^k \lambda_j y^j_d\\
    &= \bar{y}_d\text{.}
  \end{align*}
  As \(\bar{y}\) is a \(K\)-minimal point, we get with \eqref{eq:kmin}
  that \(\sum_j\lambda_j t^j\) solves the minimization problem
  \begin{equation*}
    \min\valf[\transp{\bar{w}}\datapoints]{t}\qquad\text{s.t.}\quad t\in\R\text{.}
  \end{equation*}
  According to the proposition in the beginning, the solution to this
  problem is uniquely determined.  Hence,
  \(\bar{t} = \sum_j\lambda_j t^j\).  Therefore,
  \begin{equation*}
    \transp{\bar{w}}z = \sum_{j = 1}^k \lambda_j\transp{w^j}z
    \geq \sum_{j = 1}^k \lambda_j t^j
    \geq \bar{t}
    \geq \quantile[\transp{\bar{w}}\datapoints]{\prob}\text{,}
  \end{equation*}
  where the first inequality is due to \eqref{eq:pr0}.  This proves
  \(z \in \cquantile[\transp{\bar{w}}\datapoints]{\prob}\).
\end{proof}

The classical \emph{Tukey depth regions} \cite{Tukey75} provide an
important special case of the cone quantiles $Q^-_{X,C}(p)$ and are
obtained for the choice $C=\Set{0}$, see \cite[Section
5]{HamKos18}. Since this cone has empty interior, the above results
cannot be applied directly. However, it is possible to lift the data
set as
\begin{equation*}
  \ell(X) = \Set*{
    \begin{pmatrix}
    x \\
    -e^T x
    \end{pmatrix}
    \given x \in X
  }\text{.}
\end{equation*}
The lifted data set $\ell(X)$ in $\R^{d+1}$ belongs to the hyperplane
$H=\Set{z \in \R^{d+1} \given e^T z = 0}$. For $C=\R^{d+1}_+$ we
obtain a finite representation of $Q^-_{\ell(X),C}(p)$ from a solution
$\bar T$ of \eqref{eq:vlp_d}. Defining the ``unlifted'' normals as
$\lambda(w) \coloneq (w_1-w_{d+1},\dots,w_d-w_{d+1})$, we have
$w^T \ell(X) = \lambda(w)^T X$ and we obtain a finite representation
of Tukey depth regions as
\begin{equation*}
  Q^-_{X,\{0\}}(p)
  =
  \bigcap_{(w,t) \in \bar T}
  \Set*{z \in \R^d
    \given
    \lambda(w)^T z
    \geq t
  }\text{.}
\end{equation*}

It should be noted that an alternative extension of quantiles to a
multivariate context can be found in \cite{HalPaiSim10}, where no
order~structure is used at all (e.g.,\ \(C = \Set{0}\)).  (This
article may be consulted for further references about different
approaches to multivariate extensions of the concept of quantiles.)
The \emph{quantile regions} defined there coincide with the cone
quantiles from \cite{HamKos18} in the case of the trivial ordering
cone \(C = \Set{0}\).  In \cite{HalPaiSim10}, computation of a
quantile region is carried out using parametric linear programming.
While the main result of this present note is derived directly, it
could, at least in the case \(C = \Set{0}\), also be inferred from
\cite{HalPaiSim10} as an implication of the equivalence between vector
linear programs and parametric linear programming (c.f.\
\cite{LoeWei16}).

\subsection*{Acknowledgements}
This research was motivated by a talk of Daniel Kostner and by
discussions with Andreas H. Hamel at the SKI\footnote{scientific key
  ideas} workshop in Bruneck-Brunico in March 2023.

\bibliographystyle{abbrv}
\bibliography{ref}
\end{document}